\documentstyle{article}
\input epsf
\input amssym

\def\be{\begin{eqnarray}}
\def\ee{\end{eqnarray}}

\begin{document}
\date{}
\title{\sc Flag Curvature of Invariant Randers Metrics on Homogeneous Manifolds}
\author{ E. Esrafilian\\ and\\ \vspace{1cm} H. R. Salimi Moghaddam\\
\emph{Department of Mathematics},\\ \emph{University of Isfahan}, \\ \emph{Isfahan,81746-73441-Iran.} \\
e-mail:\\
 salimi.moghaddam@gmail.com and hr.salimi@sci.ui.ac.ir}

\maketitle \abstract{} \setcounter{equation}{0} In this article we
review the recent results about the flag curvature of invariant
Randers metrics on homogeneous manifolds and by using a counter
example we show that the formula which obtained for the flag
curvature of these metrics is incorrect. Then we give an explicit
formula for the flag curvature of invariant Randers metrics on the
naturally reductive homogeneous manifolds $(G/H,g)$, where the
Randers metric induced by the invariant Riemannian metric $g$ and
an invariant vector field $\tilde{X}$ which is parallel with $g$.



\section{Introduction.}
\setcounter{equation}{0} The geometry of Finsler manifolds is one
of the interesting subjects in differential geometry which found
many physical applications (For example see \cite{[1]} and
\cite{[2]}). One of important quantities which associate to a
Finsler metric is the flag curvature which is a generalization of
the concept of sectional curvature in Riemannian geometry.\\
But, in general, the computation of the flag curvature of a
Finsler metric is very difficult. Therefore, It is very important
to find an explicit and applicable formula for the flag curvature.
In this article we want to find such explicit formula for the flag
curvature of a special type of invariant Randers metrics on
homogeneous manifolds.\\
S. Deng and Z. Hou studied invariant Finsler metrics on reductive
homogeneous manifolds and gave an algebraic description of these
metrics and obtained a necessary and sufficient condition for a
homogeneous manifold to have invariant Finsler metrics (See
\cite{[4]}). Also they studied invariant Randers metrics on
homogeneous Riemannian manifolds and used of this structure to
construct Berwald space which is neither Riemannian nor locally
Minkowskian (For more  details see \cite{[5]}). They gave a
formula for the flag curvature of invariant Randers metrics on
homogeneous manifolds in \cite{[5]}.\\
In this paper by using a counter example we show that the formula
which obtained in \cite{[5]} is incorrect. Also we explain why
this formula and also an example which gave in \cite{[4]} are
incorrect. Then we give an explicit formula for the flag curvature
of invariant Randers metrics on naturally reductive homogeneous
manifolds $(G/H,g)$, where the Randers metric induced by the
invariant Riemannian metric $g$ and an invariant vector field
$\tilde{X}$ which is parallel with $g$.

\section{Preliminaries.}
\setcounter{equation}{0}
\paragraph{2.1. Definition.}(See\cite{[6]} or \cite{[8]}) A
homogeneous space $G/H$ of a connected Lie group $G$ is called
\textit{reductive} if the following conditions are satisfied:\\

(1) In the Lie algebra $\frak{g}$ of $G$ there exists a subspace
$\frak{m}$ such that \hspace*{1.5cm}$\frak{g}=\frak{m}+\frak{h}$
(direct sum of vector subspaces).\\

(2) $ad(h)\frak{m}\subset \frak{m}$ for all $h\in H$,\\

\noindent where $\frak{h}$ is the subalgebra of $\frak{g}$
corresponding to the identity component $H_0$ of $H$ and $ad(h)$
denotes the
adjoint representation of $H$ in $\frak{g}$.\\
condition (2) implies\\

(2\'{)} $[\frak{h}, \frak{m}]\subset \frak{m}$\\

\noindent and, conversely, if $H$ is connected, then (2\'{)}
implies (2).\\
$G/H$ is reductive in either of the following cases
(See\cite{[8]}):
\begin{itemize}
    \item $H$ is compact,
    \item $H$ is connected and semi-simple.
    \item $H$ is a discrete subgroup of $G$; $\frak{h}=0$ and $\frak{m}=\frak{g}$
\end{itemize}
Let $G/H$ be a reductive homogeneous manifold with invariant
Riemannian metric $g$ which the subspace $\frak{m}$ is the
orthogonal complement  of $\frak{h}$ with respect to the inner
product on $\frak{g}$. Also let
\begin{eqnarray*}
 V=\{X\in\frak{m}|ad(h)X=X, <X,X><1, \forall h\in H\}
\end{eqnarray*}
Where $<,>$ is the inner product induced by $g$.\\
Then for any $X\in V$ there exist an invariant Randers metric on
$G/H$ by the following formula (See
\cite{[5]}):
\begin{equation}\label{F}
    F_{X}(xH,y)=\sqrt{g(xH)(y,y)}+g(xH)(\tilde{X},y), \ \ y\in
    T_{xH}(G/H)
\end{equation}
Where $\tilde{X}$ is the corresponding invariant vector field on
$G/H$ to $X$.
\paragraph{2.2. Theorem.}(See\cite{[5]}) Let $Y$ be a nonzero
vector in $\frak{m}$ and $P$ be a plane in $\frak{m}$ containing
$Y$. Then the flag curvature of the flag $(P,Y)$ in $T_0(G/H)$ is
given by \begin{eqnarray*}
  K(P,Y)=\frac{2\sqrt{g(Y,Y)}}{2\sqrt{g(Y,Y)}+g(X,Y)}K(P).
\end{eqnarray*}
Where $K(P)$ is the Riemannian curvature of $P$ of the Riemannian
metric $g$, and \begin{eqnarray*}
  K(P)=\frac{g([[Y,U]_{\frak{h}},Y],U)}{g(U,U)g(Y,Y)-g^2(U,Y)}.
\end{eqnarray*}
Where $U$ is any vector in $P$ such that $span\{Y,U\}=P$ and
$[Y,U]_{\frak{h}}$ is the orthogonal projection of $[Y,U]$ to
$\frak{h}$.\\

We will show that the above Theorem is not correct.
\section{Counter example.}
Let $G$ be a connected Lie group and $H=\{e\}$. So $G/H=G$ is a
reductive homogeneous manifold with $\frak{m}=\frak{g}$ and
$\frak{h}=0$. Also let $g$ be an invariant Riemannian metric on
$G/H=G$ and
$0=X\in \frak{m}=\frak{g}$.\\
In this case the invariant Randers metric defined by $g$ and $X$
is Riemannin because:
\begin{eqnarray*}
  F_X(xH,y)=\sqrt{g(xH)(y,y)}+g(xH)(0,y)=\sqrt{g(x)(y,y)}.
\end{eqnarray*}
Also we know that if our Finsler metric be Riemannian then the
flag curvature reduces to the familiar sectional
curvature(\cite{[3]}).\\
Under the above conditions we have
\begin{eqnarray*}
  [Y,U]_{\frak{h}}=0,
\end{eqnarray*}
because $\frak{h}=0$.\\
So we have
\begin{eqnarray*}
  g([[Y,U]_{\frak{h}},Y],U)=g([0,Y],U)=g(0,U)=0.
\end{eqnarray*}
Therefore $K(P)=0$ and so $K(P,Y)=0$.\\
By attention to this fact that in this case the flag curvature is
the same sectional curvature of $g$, we obtain the following
incorrect proposition:\\

$\clubsuit$\textit{The sectional curvature of all Lie groups with
invariant
Riemannian metrics is zero.}\\

The above proposition is false. We give the following counter
example that show there exist many Lie groups with invariant
Riemannian metrics and nonzero sectional curvatures.
\paragraph{3.1. Theorem.}(See \cite{[10]}) Let $G$ be a connected
nilpotent Lie group. Let $B$ be a positive definite symmetric
bilinear form on the Lie algebra $\frak{g}$ of $G$, and let $M$ be
the Riemannian manifold obtained by left translation of $B$ to
every tangent space of $G$. Then these are equivalent:
\begin{itemize}
    \item $M$ has a positive sectional curvature.
    \item $M$ has a negative sectional curvature.
    \item $G$ is not commutative.\hfill $\Box$
\end{itemize}
So if we consider a connected nilpotent noncommutative Lie group
$G$ with an invariant Riemannian metric, then it has nonzero
sectional curvature. (One can find another counter examples in
\cite{[7]}.)\\
By attention to the above counter example we showed that the
proposition $\clubsuit$ is incorrect so the Theorem 2.2 is
incorrect.
\section{Flag curvature of invariant Randers metrics on homogeneous manifolds.}
But why the Theorem 2.2 is incorrect? \\
If we see the proof of this Theorem we see that the authors of the
paper \cite{[5]} use of the formula
\begin{equation}\label{R}
    R(U,V)W=-[[U,V]_{\frak{h}},W]
\end{equation}
for the curvature tensor of $F$ (and $g$). But the equation
\ref{R} is a formula for tensor curvature of a special type of
invariant affine connections called the canonical affine
connection of the second kind (See Theorem 10.3 of \cite{[8]} or
Theorem 2.6 of \cite{[6]} page 193) which It is not the Riemannian
connection. So we can not use of equation \ref{R} for the
curvature tensor of the Riemannian connection.\\
Another problem of the proof of Theorem 2.2 is the formula of
$g_Y(U,V)$. The formula
$g_Y(U,V)=g(U,V)(1+\frac{g(X,Y)}{2\sqrt{g(Y,Y)}})$ which is used
in the proof of Theorem 2.2 is not correct.\\

The authors of the paper \cite{[5]} used of the equation \ref{R}
to compute the flag curvature of an example in the paper
\cite{[4]} (The example after Theorem 2.1 of the paper
\cite{[4]}). By attention to this fact that the equation \ref{R}
is not the curvature tensor of the Riemannian connection, so the
flag curvature of this example is not correct.\\

Now by attention to the above discussion, we compute the flag
curvature of invariant Randers metrics on homogeneous manifolds.
\paragraph{4.1. Definition.}(See \cite{[6]}) A homogeneous
manifold $M=G/H$ with a $G-$invariant indefinite Riemannian metric
$g$ is said to be naturally reductive if it admits an
$ad(H)$-invariant decomposition $\frak{g}=\frak{h}+\frak{m}$
satisfying the condition
\begin{eqnarray*}
 B(X,[Z,Y]_{\frak{m}})+B([Z,X]_{\frak{m}},Y)=0 \hspace{1.5cm}for  X, Y, Z \in
 \frak{m}.
\end{eqnarray*}
Where $B$ is the bilinear form on $\frak{m}$ induced by $\frak{g}$
and $[,]_{\frak{m}}$ is the projection to $\frak{m}$ with respect
to the decomposition $\frak{g}=\frak{h}+\frak{m}$.
\paragraph{4.2. Theorem.} Let $G/H$ be a homogeneous manifold with
invariant Riemannian metric $g$ and $F$ be an invariant Randers
metric defined by the $ad(H)$-invariant vector $X$,
\begin{eqnarray*}
  F(xH,Y)=\sqrt{g(xH)(Y,Y)}+g(xH)(\tilde{X},Y).
\end{eqnarray*}
Where $g(X,X)<1$, $Y\in T_{xH}G/H$ and $\tilde{X}$ is the
corresponding invariant vector field on $G/H$ to $X$. Also suppose
that the vector field $\tilde{X}$ is parallel with respect to
$g$ and $(G/H,g)$ is naturally reductive.\\
Then the flag curvature of the flag $(P,Y)$ in $T_H(G/H)$ is given
by
\begin{eqnarray*}
  K(P,Y)=\frac{A}{B-C}.
\end{eqnarray*}
Where\begin{eqnarray*}
 A&=&g(\alpha,U)+g(X,\alpha).g(X,U)-\frac{g(X,Y).g(Y,U).g(Y,\alpha)}{g(Y,Y)^{\frac{3}{2}}}\\
 &+&\frac{1}{\sqrt{g(Y,Y)}}\{g(X,\alpha).g(Y,U)+g(X,Y).g(\alpha,U)+g(X,U).g(Y,\alpha)\},
\end{eqnarray*}
\begin{eqnarray*}
  B=\{g(Y,Y)+g^2(X,Y)+2g(X,Y)\sqrt{g(Y,Y)}\}\times\\
  \{g(U,U)+g^2(X,U)-\frac{1}{\sqrt{g(Y,Y)}}\{\frac{g(X,Y).g^2(Y,U)}{g(Y,Y)}\\
  +g(X,Y).g(U,U)+2g(X,U).g(Y,U)\}\}
\end{eqnarray*}
and
\begin{eqnarray*}
  C=\{g(Y,U)(1+\frac{g(X,Y)}{\sqrt{g(Y,Y)}})+g(X,U)(g(X,Y)+\sqrt{g(Y,Y)})\}^2.
\end{eqnarray*}
Where $U$ is any vector in $P$ such that $span\{Y,U\}=P$ and in
$A$ we have
\begin{eqnarray*}
  \alpha =\frac{1}{2}[[U,Y]_{\frak{m}},Y]+[Y,[Y,U]_{\frak{h}}]
\end{eqnarray*}
($[,]_{\frak{m}}$ and $[,]_{\frak{h}}$ are the projections of
$[,]$ to $\frak{m}$ and $\frak{h}$ respectively.)\\

\noindent\textbf{\textit{Proof.}} $\tilde{X}$ is parallel with
respect to $g$ so by using Lemma 2.1 of \cite{[5]} $F$ is of
Berwald type. Also the Chern connection of $F$ and the Riemannian
connection of $g$ coincide (See \cite{[3]} page 305), therefore we
have
\begin{eqnarray*}
  R^F(U,V)W=R^g(U,V)W.
\end{eqnarray*}
Where $R^F$ and $R^g$ are the curvature tensors of $F$ and $g$
respectively. (Now let $R:=R^F=R^g.$)\\
But $g$ is naturally reductive, so by using Proposition 3.4 of
\cite{[6]} (page 202) we have:
\begin{eqnarray*}
 (R(U,V)W)_0&=&\frac{1}{4}[U,[V,W]_{\frak{m}}]_\frak{m}-\frac{1}{4}[V,[U,W]_{\frak{m}}]_{\frak{m}}\\
 &-&\frac{1}{2}[[U,V]_{\frak{m}},W]_{\frak{m}}-[[U,V]_{\frak{h}},W]
 \ \ \ for U,V,W \in \frak{m}.
\end{eqnarray*}
So we have
\begin{eqnarray*}
(R(Y,U)Y)_0&=&\frac{1}{4}[Y,[U,Y]_{\frak{m}}]_\frak{m}-\frac{1}{4}[U,[Y,Y]_{\frak{m}}]_{\frak{m}}\\
 &-&\frac{1}{2}[[Y,U]_{\frak{m}},Y]_{\frak{m}}-[[Y,U]_{\frak{h}},Y]\\
 &=&\frac{1}{2}[[U,Y]_{\frak{m}},Y]_{\frak{m}}+[Y,[Y,U]_{\frak{h}}].
\end{eqnarray*}
On the other hand we have:
\begin{equation}\label{K}
    K(P,Y)=\frac{g_Y(R_Y(U),U)}{g_Y(Y,Y).g_Y(U,U)-g_Y^2(Y,U)}.
\end{equation}
Where
\begin{eqnarray*}
  g_Y(U,V)&=&\frac{1}{2}\frac{\partial^2}{\partial s\partial
  t}\{F^2(Y+sU+tV)\}|_{s=t=0}\\
  &=&\frac{1}{2}\frac{\partial^2}{\partial s\partial
  t}\{g(Y+sU+tV,Y+sU+tV)\\
  &+&g^2(X,Y+sU+tV)\\
  &+&2\sqrt{g(Y+sU+tV,Y+sU+tV)}g(X,Y+sU+tV)\}|_{s=t=0}
\end{eqnarray*}
By a direct computation we get
\begin{eqnarray*}
  g_Y(U,V)&=&g(U,V)+g(X,U).g(X,V)-\frac{g(X,Y).g(Y,V).g(Y,U)}{g(Y,Y)^{\frac{3}{2}}}\\
  &+&\frac{1}{\sqrt{g(Y,Y)}}\{g(X,U).g(Y,V)\\
  & &+g(X,Y).g(U,V)+g(X,V).g(Y,U)\}.
\end{eqnarray*}
Therefore
\begin{eqnarray*}
  g_Y(Y,Y)&=&g(Y,Y)+g(X,Y)(g(X,Y)+2\sqrt{g(Y,Y)}),\\
  g_Y(U,U)&=&g(U,U)+g^2(X,U)-\frac{g(X,Y).g^2(Y,U)}{g(Y,Y)^{\frac{3}{2}}}\\
  &+&\frac{1}{\sqrt{g(Y,Y)}}\{g(X,Y).g(U,U)+2g(X,U).g(Y,U)\},\\
  g_Y(Y,U)&=&g(Y,U)(1+\frac{g(X,Y)}{\sqrt{g(Y,Y)}})+g(X,U)(g(X,Y)+\sqrt{g(Y,Y)}),
\end{eqnarray*}
and
\begin{eqnarray*}
  g_Y(R(Y,U)Y,U)=g_Y(\alpha,U).
\end{eqnarray*}
Where
$\alpha=\frac{1}{2}[[U,Y]_{\frak{m}},Y]_{\frak{m}}+[Y,[Y,U]_{\frak{h}}]$.\\
Combining the above formulas with the equation \ref{K} completes
the proof. \hfil$\Box$
\paragraph{4.3. Note.} In the formula obtained in Theorem 4.2 if
we assume that $U$ is orthogonal to $Y$ with respect to $g$ then
we can obtain a simpler formula.
\paragraph{4.4. Corollary.} In a special case we can assume a
Lie group $G$ as a reductive homogeneous space with $H=\{e\}$,
$\frak{h}=\{0\}$ and $\frak{m}=\frak{g}$, then our formula for the
flag curvature is simpler because in this case we have
\begin{eqnarray*}
 \alpha=\frac{1}{2}[[U,Y],Y].
\end{eqnarray*}

\end{document}